\documentclass[paper=letter,11pt]{scrartcl}

\usepackage[utf8]{inputenc} \usepackage[T1]{fontenc}    
\usepackage{natbib}
\setcitestyle{numbers, compact, open={[},close={]}}

\usepackage{amssymb,amsmath,amsthm,sectsty,url,graphicx}
\usepackage{hyperref}       \usepackage{lmodern}
\usepackage{microtype}
\usepackage{fullpage}

\usepackage{latexsym}
\usepackage{graphicx}
\usepackage{mathptmx}

\usepackage{subcaption}		
\usepackage{amsmath}
\usepackage{amsfonts}
\usepackage{amssymb}
\usepackage{amsbsy}
\usepackage{amsthm}
\usepackage{float}

\usepackage{xspace}
\newcommand\SPSA{\texttt{SPSA}\xspace}
\newcommand\STR{\texttt{STRONG}\xspace}
\newcommand\STRONE{\texttt{STRONGstg1}\xspace}
\newcommand\GS{\texttt{GradSearch}\xspace}
\newcommand\RS{\texttt{RandomSearch}\xspace}
\newcommand\NM{\texttt{Nelder-Mead}\xspace}
\newcommand\SIMOPT{\texttt{SimOpt}\xspace}
\begin{document}
	
	\title{Comparing the Finite-Time Performance of Simulation-Optimization Algorithms}
	
		\author{		Naijia (Anna) Dong \and
		David J. Eckman \and
		Matthias Poloczek \and
		Xueqi Zhao \and
		Shane G. Henderson\thanks{School of Operations Research and Information Engineering, Cornell University, 206 Rhodes Hall, Ithaca, NY 14853, USA}
	}
	\date{}
	
	\maketitle
	
	\section*{ABSTRACT}
	We empirically evaluate the finite-time performance of several simulation-optimization algorithms on a testbed of problems with the goal of motivating further development of algorithms with strong finite-time performance.
We investigate if the observed performance of the algorithms can be explained by properties of the problems, e.g., the number of decision variables, the topology of the objective function, or the magnitude of the simulation error. 	
	\section{INTRODUCTION}
	The practice of simulation optimization (SO) deals with optimizing a real-valued objective function that cannot be evaluated exactly, but instead must be estimated via simulation.
In addition to the challenge of estimating the objective function, structural properties of the objective function, e.g., continuity, differentiability, or convexity, may be unknown.
On account of this, SO algorithms are often designed to solve a broad class of problems without exploiting any structure of the objective function.

In the SO literature, one commonly sees theoretical results on the asymptotic performance of an algorithm.
It is often shown that an algorithm will converge to a local or global optimizer as the simulation effort approaches infinity.
Some results further specify a rate at which an algorithm converges once within a neighborhood of an optimizer.
Unfortunately, the asymptotic regimes in which these results hold likely require an amount of computational effort that exceeds practical budgets, making the asymptotic results less useful to practitioners.
For a practitioner deciding which algorithm to use for a particular problem, understanding the finite-time performance of an algorithm is more meaningful.

The SO community lags behind other research communities when it comes to having an established testbed of problems and developing metrics for comparing the finite-time performance of algorithms \cite{pasupathy:06}.
Moreover, a comprehensive comparison of SO algorithms has not been done on a large testbed \cite{amaran:16}.
As a first step toward such a comparison, we implement several popular SO algorithms and test them on a subset of problems from \href{http://simopt.org}{\SIMOPT} (\citet{henderson:11}), a growing library of SO problems developed by \citet{pasupathy:11}.  We evaluate the finite-time performance of the algorithms and discuss insights into the types of problems on which the algorithms might be expected to perform well.
An objective of our study is to spur the development of SO algorithms that have strong finite-time performance for various classes of problems.
We also hope to encourage further contributions to the \href{http://simopt.org}{\SIMOPT} library so that a more comprehensive testbed of SO problems is available to researchers for evaluating new algorithms.

	\section{EVALUATING FINITE-TIME PERFORMANCE}\label{sec:finite_perf}
	In deterministic optimization, an algorithm's performance on a given problem is usually measured by either the number of function evaluations or the wall-clock time needed to find the optimal solution or to get within a specified tolerance.
In this way, one can easily compare algorithms.
Applying this approach to SO algorithms runs into several challenges.
Firstly, the optimal solutions to SO problems are often unknown or have no certificate of optimality. Secondly, sampling error makes it harder to determine whether the objective function value of a solution is within a given tolerance of the optimal objective function value.
Instead, SO algorithms can be more fairly compared by fixing a given simulation budget---either the wall-clock time or the number of objective function evaluations---and evaluating the objective function at the estimated best solution visited within the budget.

Although measuring a computational budget in terms of wall-clock time may make sense from a practical standpoint, the resulting performances are platform dependent.
On the other hand, using the number of objective function evaluations as a measure of time has its own issues; see also~\citet{pasupathy:06}.
Firstly, a potentially significant portion of computation effort may go unmeasured if, for example, the constraints are stochastic or gradient information is calculated without taking additional objective function valuations.
Secondly, for steady-state simulations that involve simulating a single, long sample path, counting one replication as one objective function evaluation may be misleading.

In our experiments, we specify the simulation budget in terms of the number of objective function evaluations.
All of the problems we study have deterministic constraints and the simulations all have finite horizons.
Moreover, we count all objective function evaluations towards the budget, including function evaluations used to obtain gradient estimates, e.g., via finite differences; see also Sect.~\ref{sect_algos}.

For one macroreplication of an algorithm, let $Z(n)$ denote the \emph{true} objective function value of the \emph{estimated} best solution visited in the first $n$ objective function evaluations.
Because the estimated best solution $X(n)$ is random, $Z(n)$ is a random variable.
Conditional on the solution $X(n)$, the objective value $Z(n)$ is not random, but we probably cannot compute it exactly because we use simulation to evaluate the objective function.
In our experiments, we get fairly precise estimates of $Z(n)$---conditional on $X(n)$---by running additional simulations in a post-processing step.
These replications are not counted towards the algorithm's budget.
Plotting $Z(n)$ as a function of $n$ shows how the objective value of the estimated best solution changes as the algorithm progresses \cite{pasupathy:06}.
The $Z(n)$ curve has upward or downward jumps whenever a new solution is discovered that is believed to be the best.

Plotting the $Z(n)$ curve for a single macroreplication has limited value since the location of the curve is itself random.
Instead, it is more informative to run several macroreplications of an algorithm and aggregate the results.
In our experiments, we plot the mean performance curve $\bar{Z}(n) := m^{-1}\sum_{i=1}^m Z_i(n)$ for all $n$ where $Z_i(n)$ is the performance associated with the $i$th of $m$ macroreplications.
An argument can also be made for considering the \emph{median} performance of an algorithm as it is less sensitive to outliers in performance.
The empirical cumulative distribution function (cdf) of~$Z(n)$ also contains a great deal of information about how an algorithm performs, including its variability.
Although plotting the empirical cdf of $Z(n)$ for a single algorithm on a single problem is straightforward, comparing empirical cdfs for multiple algorithms is challenging.

By testing algorithms on only a modest number of SO problems, we were able to present in this paper many of the plots of $\bar{Z}(n)$.
Ideally, we would like to graphically compare the performances of algorithms \emph{across} a large number of problems as is done in the deterministic optimization community using performance profiles \cite{dolan:02}.
Adapting performance profiles to simulation optimization remains an active area of research with a notable challenge being how exactly to define the performance ratio.
We believe that performance profiles have great potential for future comparisons of SO algorithms, especially when a large testbed of problems is available.

 	
	\section{ALGORITHMS AND PROBLEMS}

\subsection{Algorithms}
\label{sect_algos}

We suppose that when evaluating a solution $x$, an algorithm observes~$y(x) = f(x) + \varepsilon(x)$ where~$f(x)$ is the value of the objective function~$f$ at~$x$ and~$\varepsilon(x)$ is the observational noise associated with simulating~$x$.
We further assume that $\varepsilon(x)$ has a mean of zero and finite variance.
In our implementations, each of the algorithms take 30 samples of a given solution to estimate its objective value.

All problems specify a deterministic domain~$\mathcal{D} \subseteq \mathbb{R}^d$ where~$d$ is the number of decision variables. 
Unless otherwise stated, the starting solution~$x_1$ is drawn from within~${\cal D}$ according to some probability distribution.
When the domain was bounded, we used a uniform distribution over~${\cal D}$ and when the domain was unbounded, we used either independent exponential or Laplace distributions for generating the components of $x_1$.
In the case that an algorithm decides to simulate a solution~$x_{k+1} \notin {\cal D}$, we move it to the boundary of~${\cal D}$ on the line connecting~$x_{k+1}$ and the previous solution~$x_k$.

For our experiments, we selected a variety of well-known SO algorithms and compared them with two baseline methods: 
random search and gradient-based search.
We provide a high-level description of the algorithms here and refer the interested reader to our \href{https://bitbucket.org/poloczek/finitetimesimopt}{public repository} (\citet{poloczek:17}) for the codes.

\subsection*{Random Search}
Random search (\RS) iteratively evaluates solutions that are drawn from~$\cal D$ according to some fixed probability distribution until the simulation budget is exhausted.

\subsection*{Gradient Search with Random Restarts}
The gradient-based search algorithm (\GS) approximates the gradient of the objective function at a given solution $x_k$ via~$d$ finite differences. That is, the~$j$th component of the approximation $\hat{g}(x_k)$ is $$\hat{g}_j(x_k) = \frac{y(x_{k} + h_j e_j) - y(x_{k} - h_j e_j)}{2 h_j},$$
where for iterations $k \ge 2$,~$h_j:=\sqrt{\mathrm{Var}(y(x_1))}/(\sqrt{2r} \cdot \hat{g}_j(x_1))$, $\mathrm{Var}(y(x_1))$ is the \emph{estimated} variance of the noise $\varepsilon(x_1)$, $r$ is the number of replications taken at~$x_1$, and~$e_j$ is the~$d$-dimensional unit vector whose~$j$th entry is one. 
In the first iteration, we randomly draw two solutions~$x_1$ and~$x'$ from $\mathcal{D}$ and set~$h_j = \min_{i} \|x_1[i] - x'[i]\|/3$ for all~$j \in \{1,\ldots,d\}$, in order to compute~$\hat{g}_j(x_1)$. 
In a given iteration $k$, we first test a step size of~$c_k = 2$ by evaluating $y(x_k + \hat{g}_k \cdot c_k)$ (for maximization problems).
If this yields a better solution than $x_k$, we set $x_{k+1} = x_k + \hat{g}_k \cdot c_k$.
Otherwise we iteratively divide $c_k$ by~$2$ and test again, until either a better solution is found or $c_k$ is too small, in which case we choose~$x_{k+1}$ to be a random point.

In order to prevent \GS from becoming trapped at a local optimum, the algorithm restarts from a randomly chosen solution if \emph{all} of the following conditions are met:
\begin{itemize}
\item[i)] $y(x_k) - y(x_{k+1}) < \tau \cdot (1 + \|y(x_{k+1})\|)$, \hspace{2mm} ii) $\quad\|x_k - x_{k+1}\| < \sqrt{\tau} \cdot (1 + \|x_{k+1}\|)$,
\item[iii)] $\|\hat{g}_{k}\| < \sqrt[3]{\tau} \cdot (1 + \|y(x_{k+1})\|)$, and \hspace{2mm} iv)  $\quad\|\hat{g}_{k}\| < \mathrm{Var}(y(x_k))$, 
\end{itemize}
where~$\tau$ is a constant chosen by the user. In our experiments, we set~$\tau = 10^{-4}$.

\subsection*{The Simultaneous Perturbation Method}
The Simultaneous Perturbation Method (\SPSA)~\cite{spall:92,spall:98} also performs an iterative search that starts from a randomly selected solution~$x_1$ and approximates the gradient in each step.
However, instead of using the method of finite differences, \SPSA relies on a simultaneous perturbation approximation of the gradient: the $j$th component of $g(x_k)$ is approximated as
$$\hat{g}_j(x_k) = \frac{y(x_{k} + c_k \Delta_{k}) - y(x_{k} - c_k \Delta_{k})}{2 c_k \Delta_{k}(j)},$$
where~$\Delta_{k}$ is a carefully chosen $d$-dimensional random vector; see \citet{spall:98} for details.
\SPSA requires only two function evaluations per iteration to approximate the gradient, whereas~$\GS$ requires~$2d$.
For the step-length sequence, we used the ``automatic gain selection'' implementation of \citet{spall:01}.

\subsection*{The Stochastic Trust-Region Response-Surface Method}
The Stochastic Trust-Region Response-Surface Method (\STR) of~\citet{chang:13} approximates the unknown objective function by a series of local models, where each model is sufficiently accurate within its corresponding \emph{trust region}.
In iteration~$k$, \STR constructs a local model~$r_k$ that is believed to resemble the objective function within a trust region ${\cal B}(x_k,\Delta_k)$ of radius $\Delta_k$ centered at $x_k$.
It then computes a point $x^\ast_k \in {\cal B}(x_k,\Delta_k)$ that is ``close'' to the optimum of $r_k$ over ${\cal B}(x_k,\Delta_k)$.
The point~$x^\ast_k$ is accepted and becomes the center point~$x_{k+1}$ for the next iteration if the following two tests are both passed:
\begin{enumerate}
\item The improvement in objective value~$y(x_k) - y(x^\ast_k)$ (for minimization problems and estimated by sampling) is sufficiently large compared to the predicted improvement~$r_k(x_k) - r_k(x^\ast_k)$.\item The improvement in objective value is statistically significant, accounting for observational noise.
\end{enumerate}
If the local model fits the observed data at~$x_k$ and~$x^\ast_k$ well---indicated by satisfactory results to the above tests---then the radius of the trust region for the next iteration is increased.
Otherwise it stays the same or, in the case of a poor prediction, is decreased.

Our implementation of~\STR follows~\citet{chang:13}, with the exception that we do not apply design of experiments to select the evaluation points for fitting.
Instead, we use central finite differences to estimate the gradient and a BFGS update to estimate the Hessian.
We also tested a version (\STRONE) that does not fit second-order models, i.e.,\ it only applies the first stage of \STR.

\subsection*{The Nelder-Mead Algorithm}
The algorithm of~\citet{nm65} (\NM) iteratively maintains a simplex of~$d+1$ vertices whose centroid is denoted by $x_k$. In iteration~$k$, the vertex with the worst observed objective value, say~$z_k$, is reflected through the centroid of the remaining~$d$ vertices to obtain a new point~$z_k'$ that is then sampled.
If the observed value~$y(z_k')$ is worse than the values previously observed at the remaining~$d$ vertices of the simplex, the simplex is contracted and~$z_k'$ is chosen closer to the centroid.
Otherwise~$z_k$ is removed from the simplex and~$z_k'$ is added as $(d+1)$st vertex.
We implemented the improvements suggested by~\citet{barton:96} for accelerating the convergence of the algorithm when function evaluations are noisy.

\subsection{Benchmark Problems}
\label{section_benchmarks}
We have compiled a testbed of 15 problems: 12 were taken from \href{http://simopt.org}{\SIMOPT} and the other three (POMDPController, RoutePrices, and TollNetwork) were developed during the course of this project and will be made available on \href{http://simopt.org}{\SIMOPT} upon publication. 
The characteristics of the problems are summarized in Table~\ref{table_problems}. 
\begin{table}[h]
  \centering
  \caption{The benchmark problems and their characteristics: The first column gives the abbreviations used in this paper and hyperlinks to their respective entries on \href{http://simopt.org}{\SIMOPT}.} 
    \label{table_problems}
    \begin{tabular}{llcc}
    \textbf{Problems} & \textbf{Name on SimOpt} & \multicolumn{1}{c}{\textbf{Dimension}} & \textbf{Optimal Solution} \\
    \hline
    \href{http://simopt.org/wiki/index.php?title=Ambulances_in_a_square}{Ambulance} & Ambulances in a Square & {6} & Unknown \\
    \href{http://simopt.org/wiki/index.php?title=Continuous_Newsvendor}{CtsNews} & Continuous Newsvendor & {1} & Known \\
    \href{http://simopt.org/wiki/index.php?title=Dual_Sourcing}{DualSourcing} & Dual Sourcing & {2} & Unknown \\
    \href{http://simopt.org/wiki/index.php?title=Economic-Order-Quantity_Model}{EOQ}   & Economic-Order-Quantity & {1} & Known \\
    \href{http://simopt.org/wiki/index.php?title=Facility_Location}{FacilityLocation} & Facility Location & {4} & Unknown \\
    \href{http://simopt.org/wiki/index.php?title=M/M/1_Metamodel}{MM1} & M/M/1 Metamodel & {3} & Known \\
    \href{http://simopt.org/wiki/index.php?title=A_Multimodal_Function}{MultiModal} & A Multimodal Function & {2} & Known \\
    \href{http://simopt.org/wiki/index.php?title=Parameter_Estimation:_2D_Gamma}{ParameterEstimation} & Parameter Estimation: 2D Gamma & {2} & Known \\
    \href{https://people.orie.cornell.edu/dje88/POMDP.pdf}{POMDPController} & Optimal Controller for a POMDP & {10} & Unknown \\
        \href{http://simopt.org/wiki/index.php?title=Optimization_of_a_Production_Line}{ProductionLine} & Optimization of a Production Line & {3} & Unknown \\
    \href{http://simopt.org/wiki/index.php?title=GI/G/1_Queue}{QueueGG1} & GI/G/1 Queue & {1} & Unknown \\
    \href{http://simopt.org/wiki/index.php?title=Rosenbrock\%27s_Function}{Rosenbrock} & Rosenbrock's Function & {40} & Known \\
    \href{https://people.orie.cornell.edu/dje88/RoutePrices.pdf}{RoutePrices} & Route Prices for Mobility-on-Demand & {12} & Unknown \\
        \href{http://simopt.org/wiki/index.php?title=SAN+Duration}{SAN}   & SAN Duration & {13} & Unknown \\
    \href{https://people.orie.cornell.edu/dje88/TollNetwork.pdf}{TollNetwork} & Toll Road Improvements& {12} & Unknown \\
        \hline
    \end{tabular}  \label{tab:addlabel}\end{table} 
 	
	\section{RESULTS AND DISCUSSION}
	For every problem, we ran 30 macroreplications of each algorithm.
Each macroreplication produced a sequence of estimated best solutions
$X(n)$, where $n$ ranges over the replication budget as specified by
the problem's description on \href{http://simopt.org}{\SIMOPT}.
As a post-processing step, we averaged 30 replications at each solution
$X(n)$ to obtain estimates of the objective function $Z(n)$. These
post-processing replications were independent of those used to identify the sequence of
solutions $X(\cdot)$, and they use common random numbers across all algorithms.
We then averaged the 30 estimates of $Z(n)$ to produce the $\bar{Z}(n)$ curve. Since \SPSA uses the replication budget $N$ as an input, we reran the algorithm with different values of $N$ to produce each point of the $\bar{Z}(n)$ curve.
We also calculated 95\% normal confidence intervals around $\bar{Z}(n)$.
In some instances, we show plots of the median performance and the
first and third quartiles of the 30 (macroreplication) samples of $Z(n)$.
We organize our discussion into groupings of problems with similar patterns in their $\bar{Z}(n)$ plots.

\subsection*{CtsNews, MM1, ParameterEstimation, and QueueGG1}
\begin{figure} [h]
	\centering
	\begin{subfigure}[b]{0.48\textwidth}
		\centering
		\includegraphics[width=\textwidth]{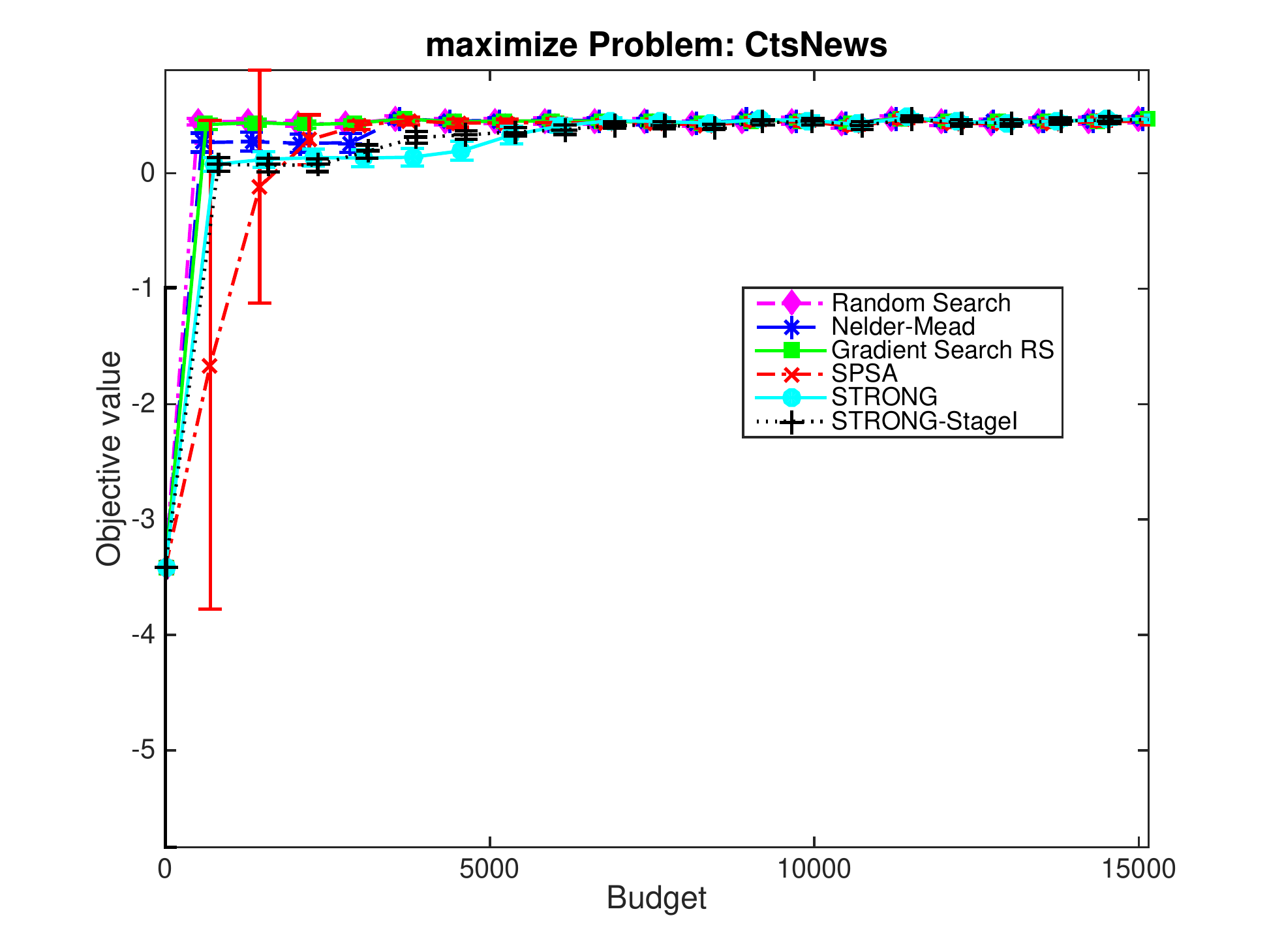}
					\end{subfigure}
	\begin{subfigure}[b]{0.48\textwidth}
		\centering
		\includegraphics[width=\textwidth]{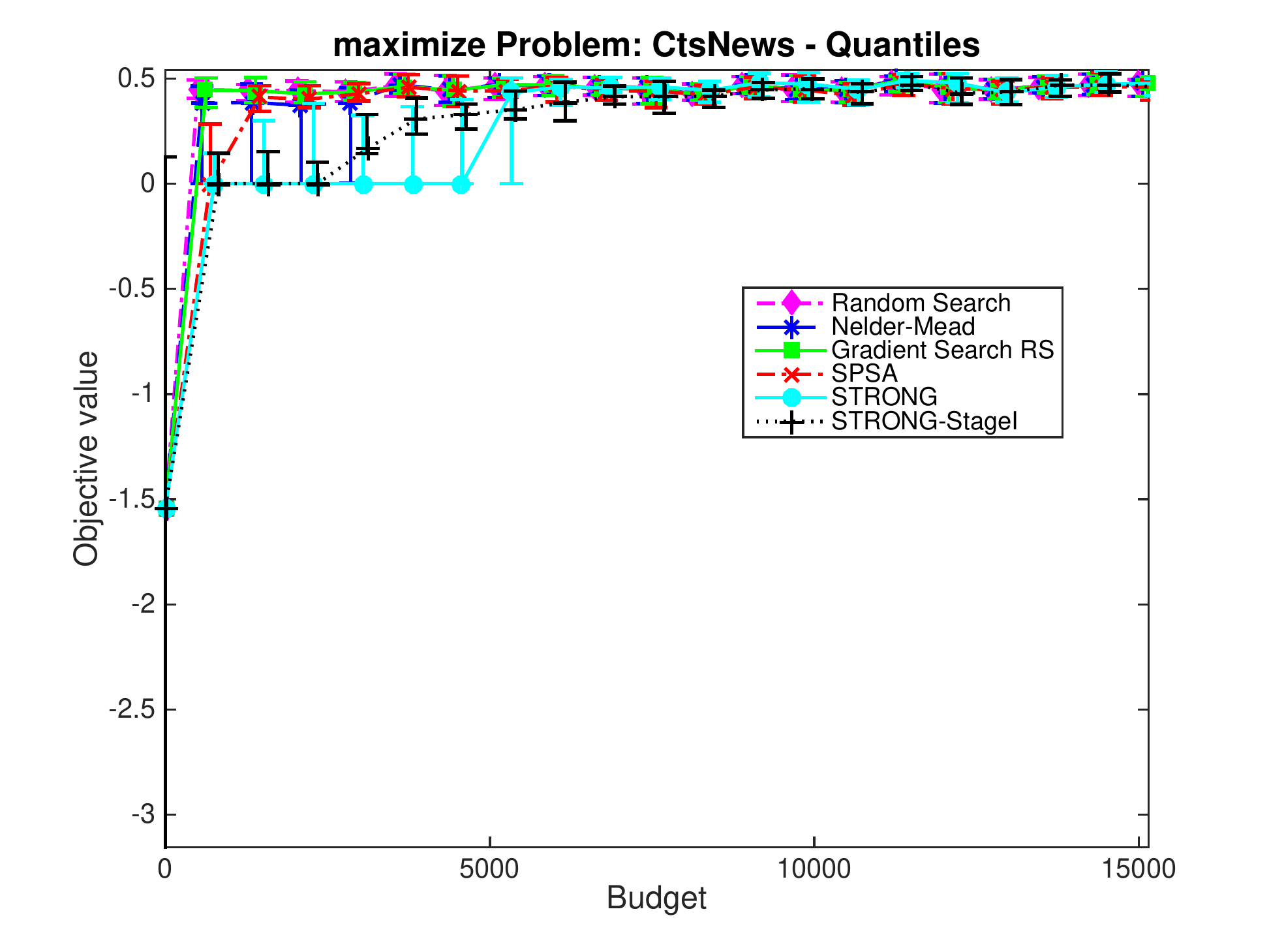}
					\end{subfigure}
	\caption{CtsNews: (l) Average Performance, (r) Performance Quantiles}
	\label{fig:CtsNews}
	\label{fig:CtsNewsQ}
\end{figure}

All algorithms work well on these low-dimensional problems, quickly converging to good solutions as illustrated in Figure~\ref{fig:CtsNews} for the problem CtsNews.
The plotted quantiles show that the performances of \NM and \STRONE are initially more variable than is suggested by the confidence intervals.
We also see that \SPSA's performance has a large variance for small
budgets, as reflected in the width of the confidence intervals. The
high variance is not seen in the corresponding quantile plot,
suggesting that \SPSA has occasional very poor performance that shows
up only in the below-0.25 quantiles (the worst quarter of the macroreplications).

\begin{figure} [h]
	\centering
	\begin{subfigure}[b]{0.48\textwidth}
		\includegraphics[width=\textwidth]{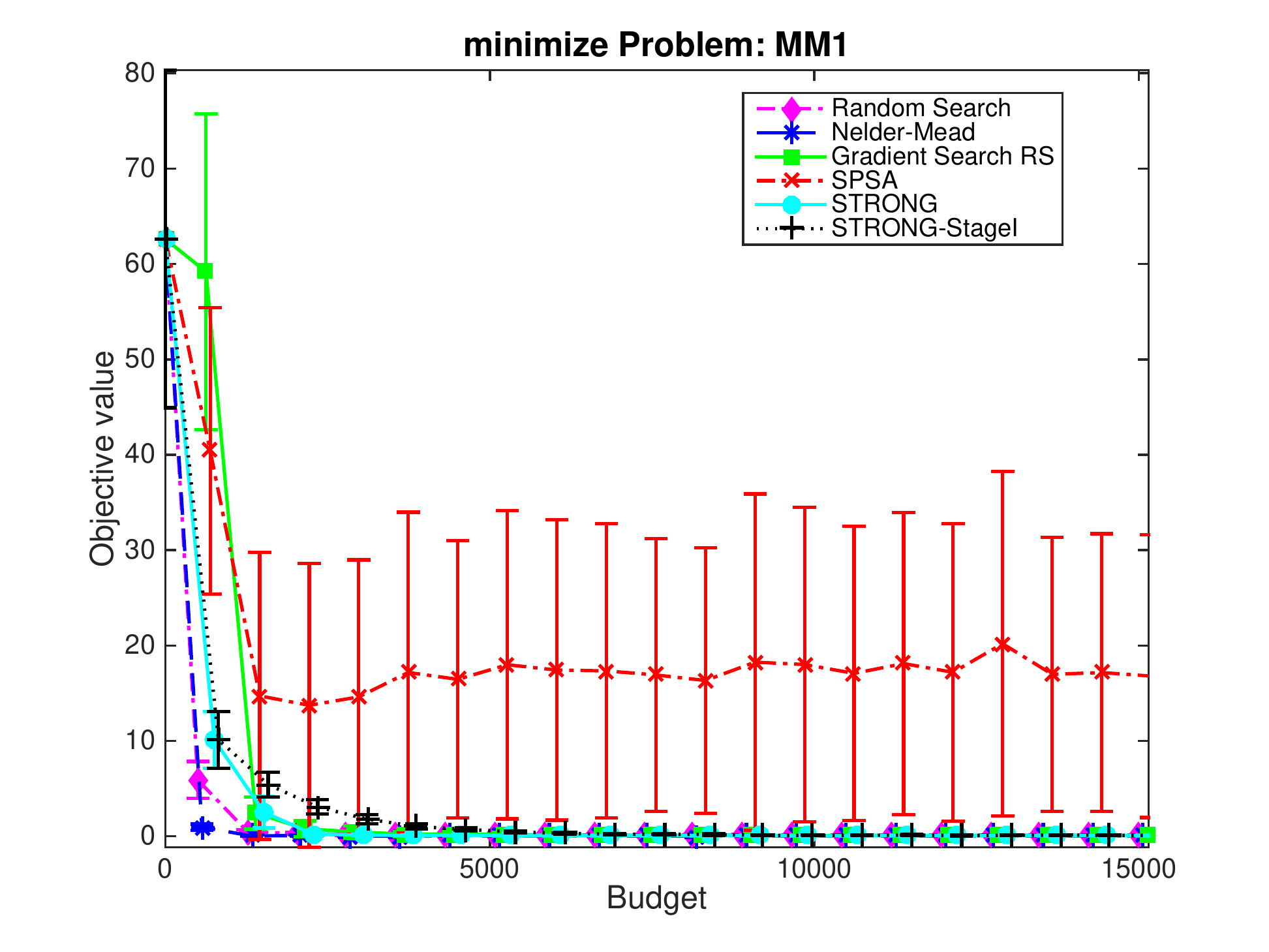}
					\end{subfigure}
	\begin{subfigure}[b]{0.48\textwidth}
		\includegraphics[width=\textwidth]{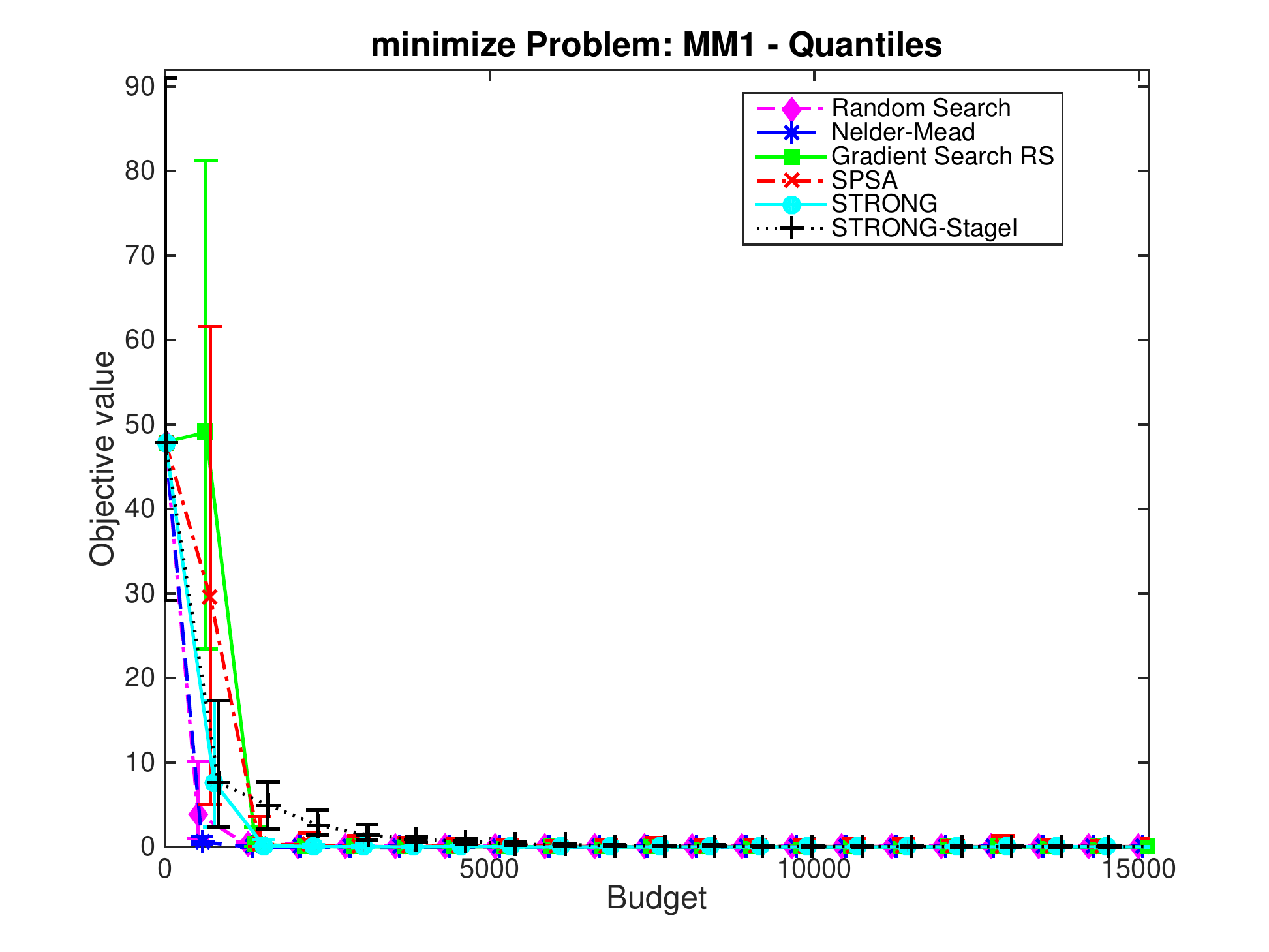}
					\end{subfigure}
	\caption{MM1: (l) Average Performance, (r) Performance Quantiles}
	\label{fig:MM1Q}\label{fig:MM1}
\end{figure}

For the problem MM1, we noticed that the performance of \SPSA was highly variable and worse than those of the other algorithms, as shown in Figure~\ref{fig:MM1}.
However, the median performance of \SPSA is competitive with those of the other algorithms.
\SPSA may be sensitive to the initial solution and occasionally can fail to make progress.
This aspect of \SPSA's performance appeared in other problems; see the discussion of POMDPController.

\subsection*{EOQ, DualSourcing}
\begin{figure} [H]
	\centering
	
	\begin{subfigure}[b]{0.48\textwidth}
		\includegraphics[width=\textwidth]{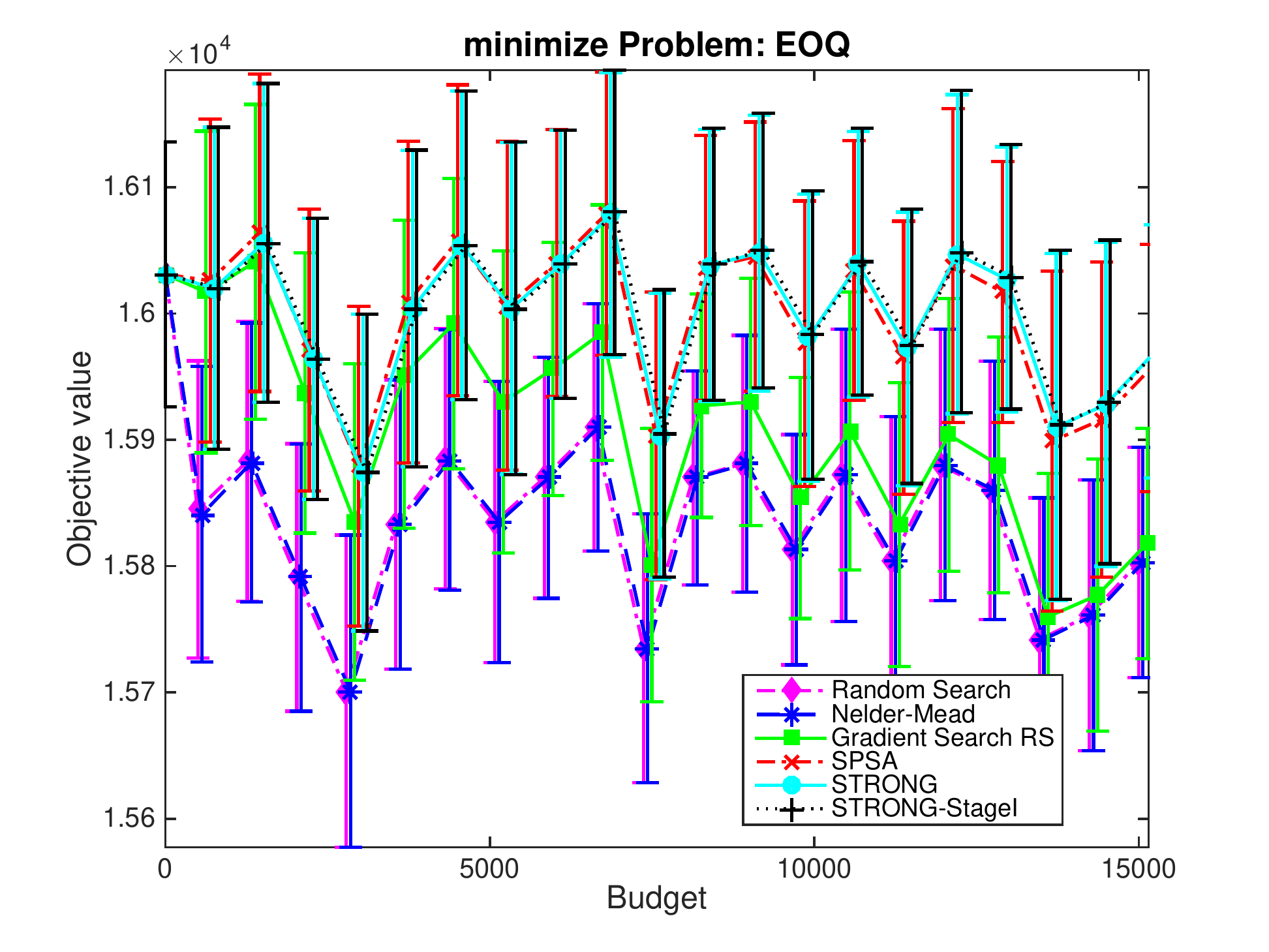}
					\end{subfigure}
	\begin{subfigure}[b]{0.48\textwidth}	
		\includegraphics[width=\textwidth]{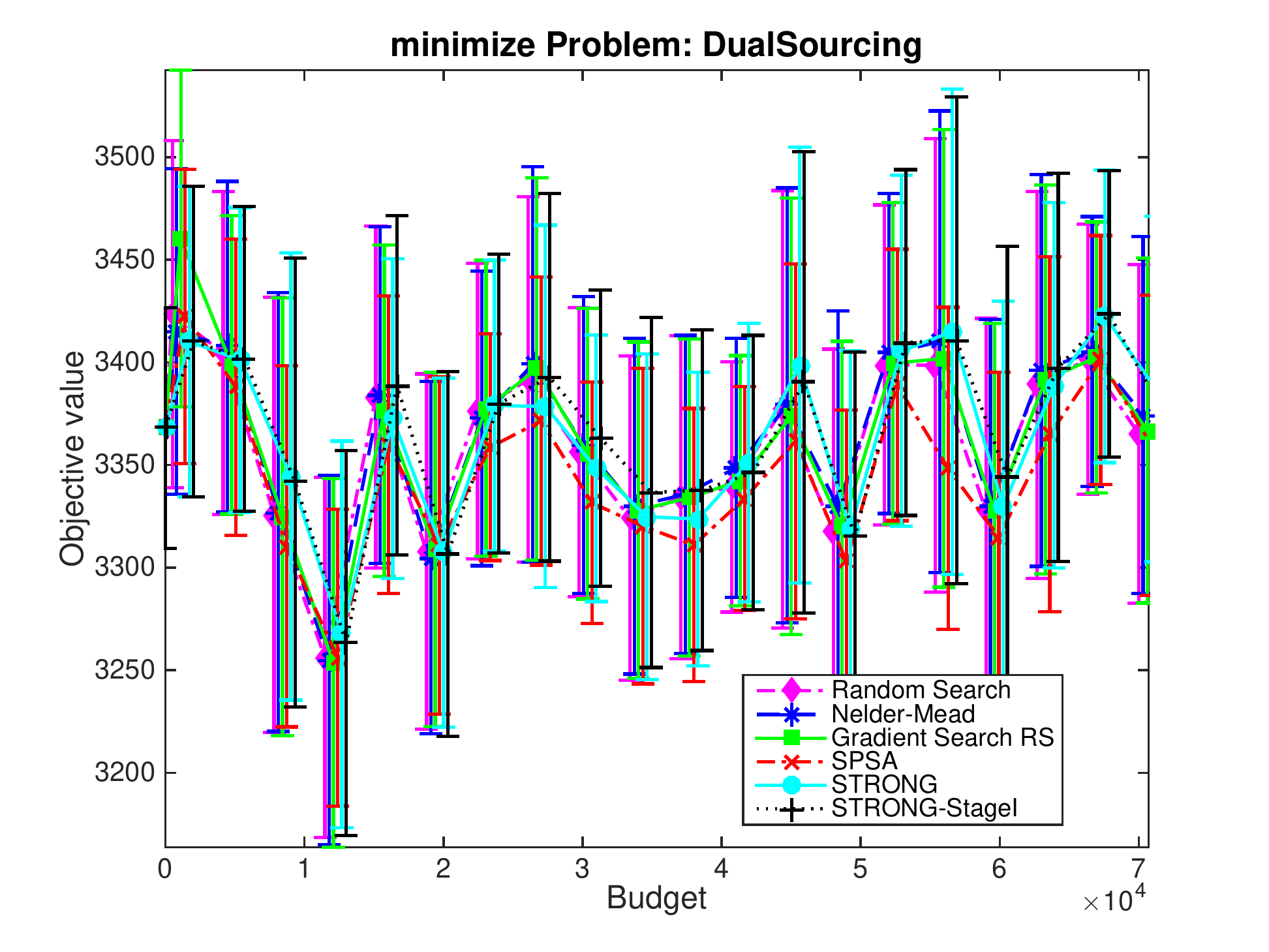}
					\end{subfigure}
	\caption{ (l) EOQ: Average Performance, (r) DualSourcing: Average Performance}
	\label{fig:EOQ}
	\label{fig:DualSourcing}
\end{figure}

\begin{figure} [H]
	\centering
	
	\begin{subfigure}[b]{0.48\textwidth}
		\includegraphics[width=\textwidth]{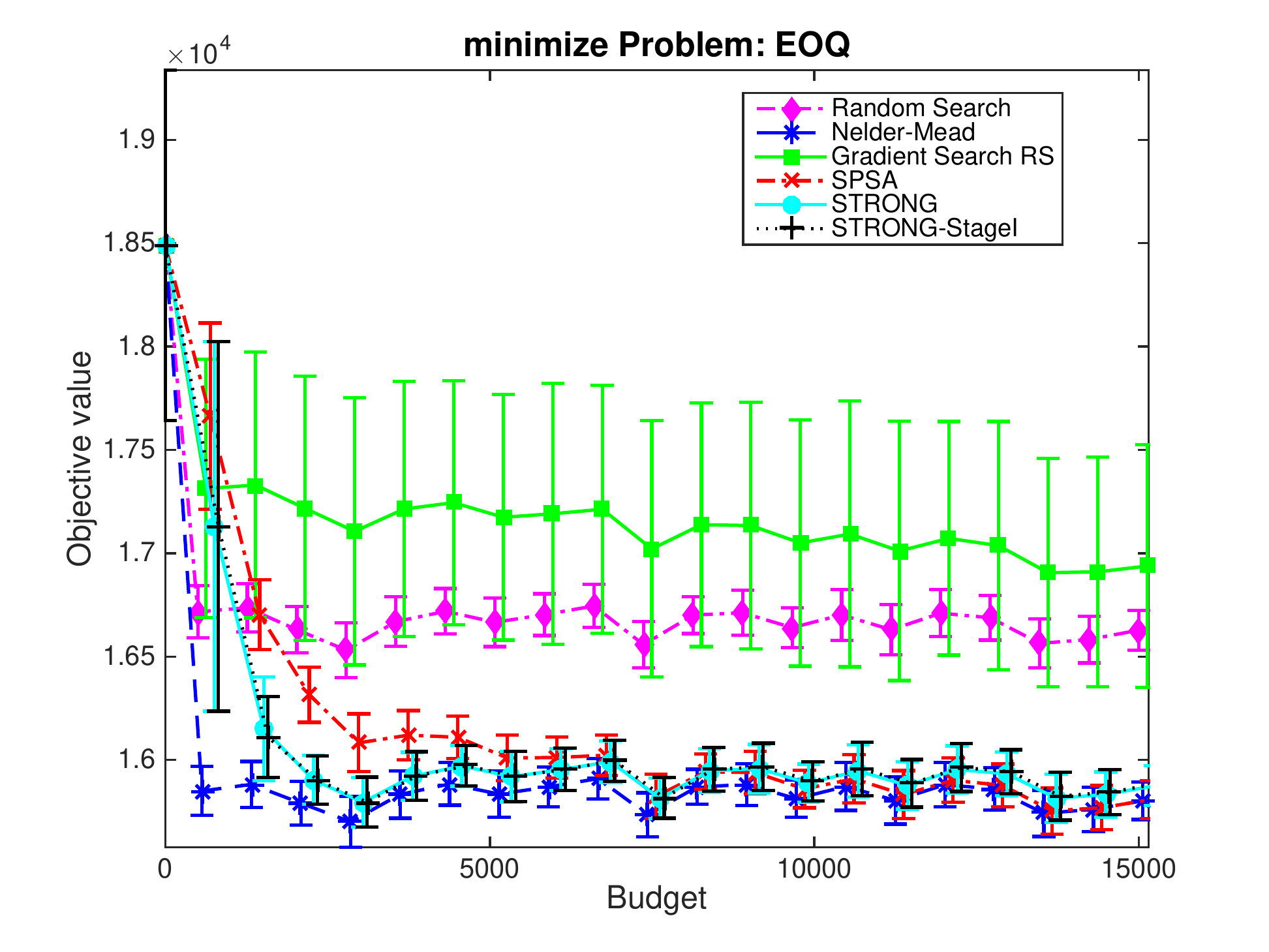}
					\end{subfigure}
	\begin{subfigure}[b]{0.48\textwidth}	
		\includegraphics[width=\textwidth]{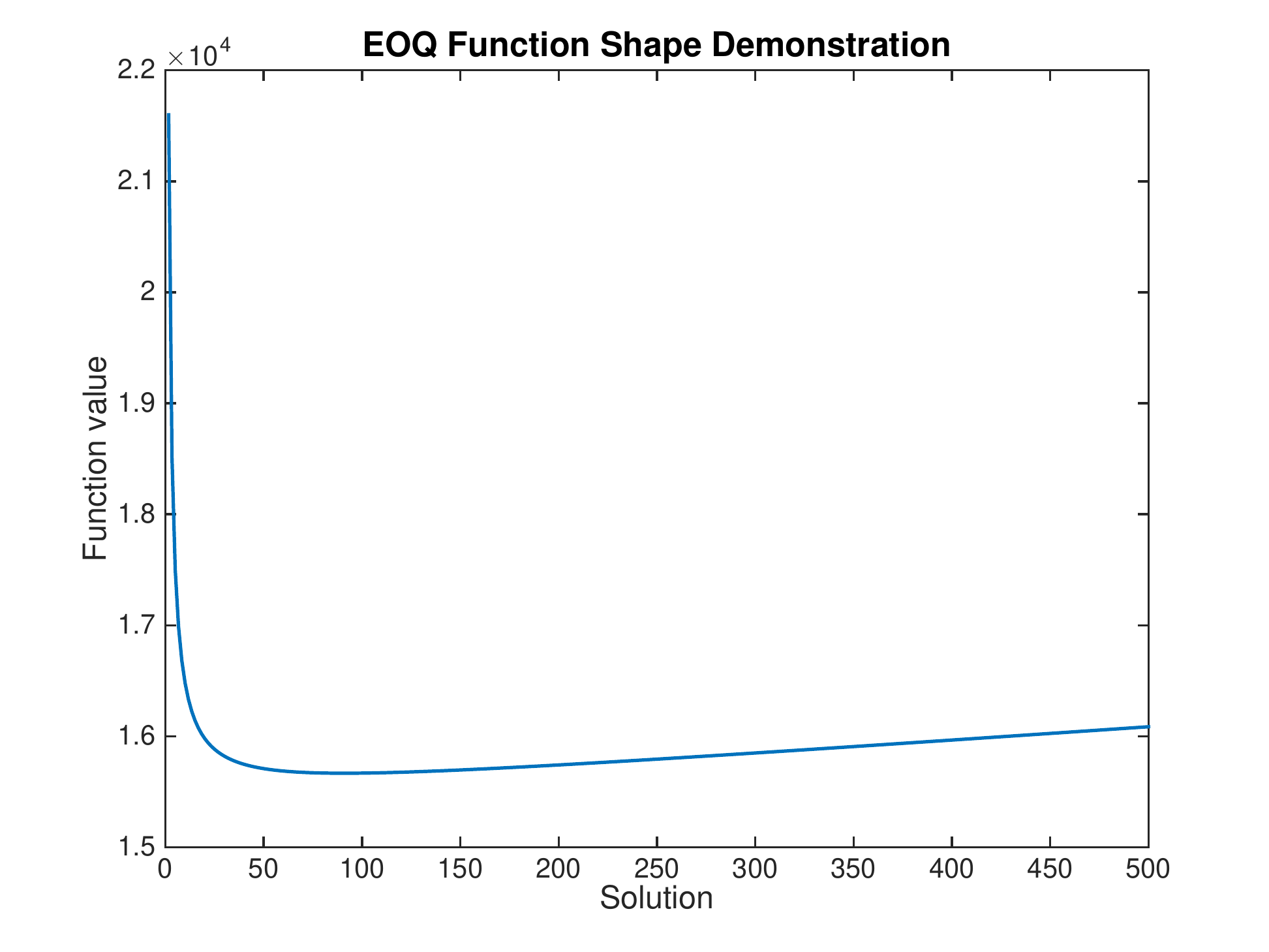}
					\end{subfigure}
	\caption{EOQ: (l) Average Performance w/ Bad Initial Solutions, (r) The Objective function}
	\label{fig:EOQ10}
	\label{fig:EOQ-Function}
\end{figure}

\begin{figure}[h]
	\centering
		\begin{subfigure}[b]{0.48\textwidth}
		\includegraphics[width=\textwidth]{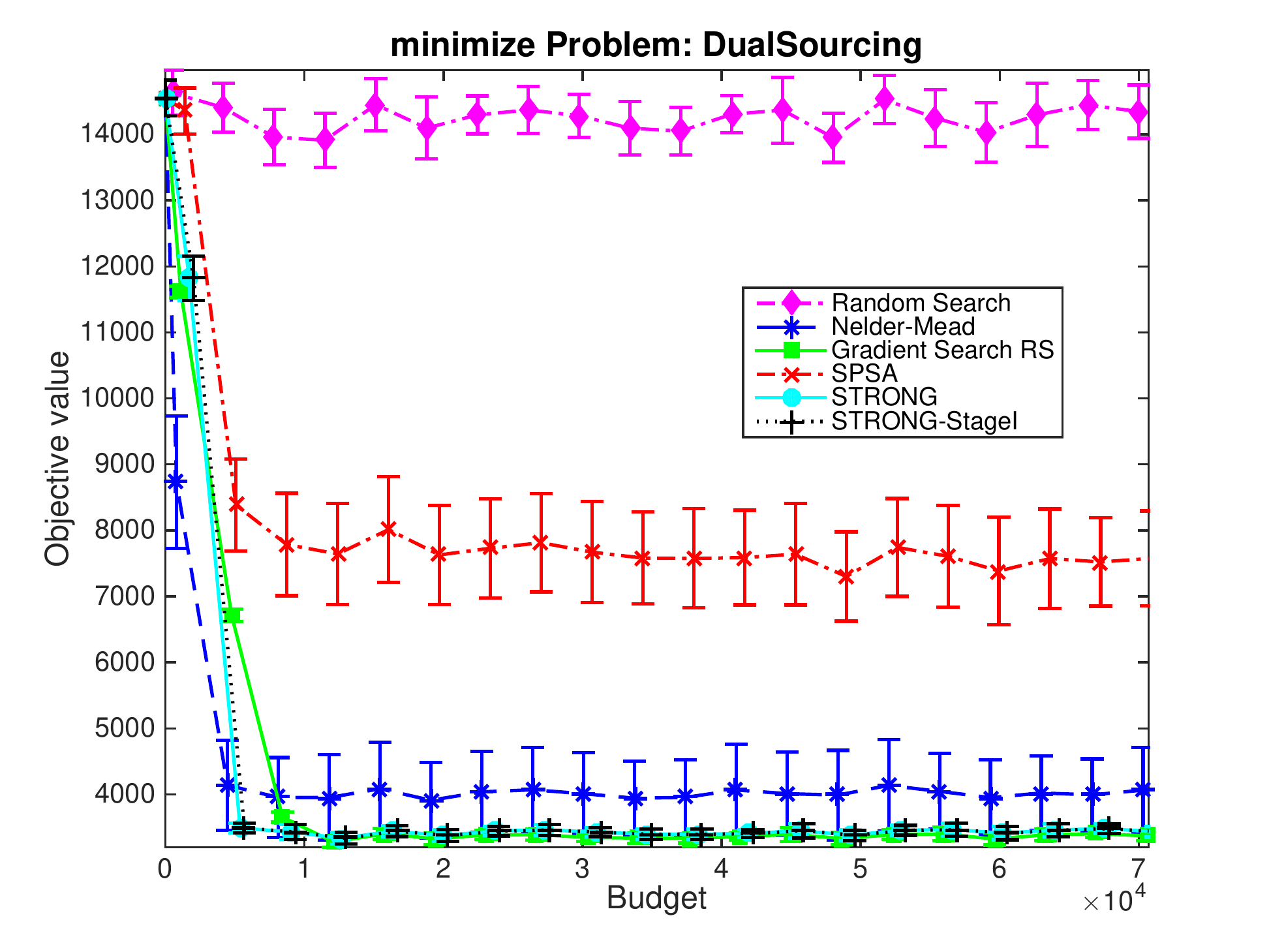}
					\end{subfigure}
	\begin{subfigure}[b]{0.48\textwidth}
		\includegraphics[width=\textwidth]{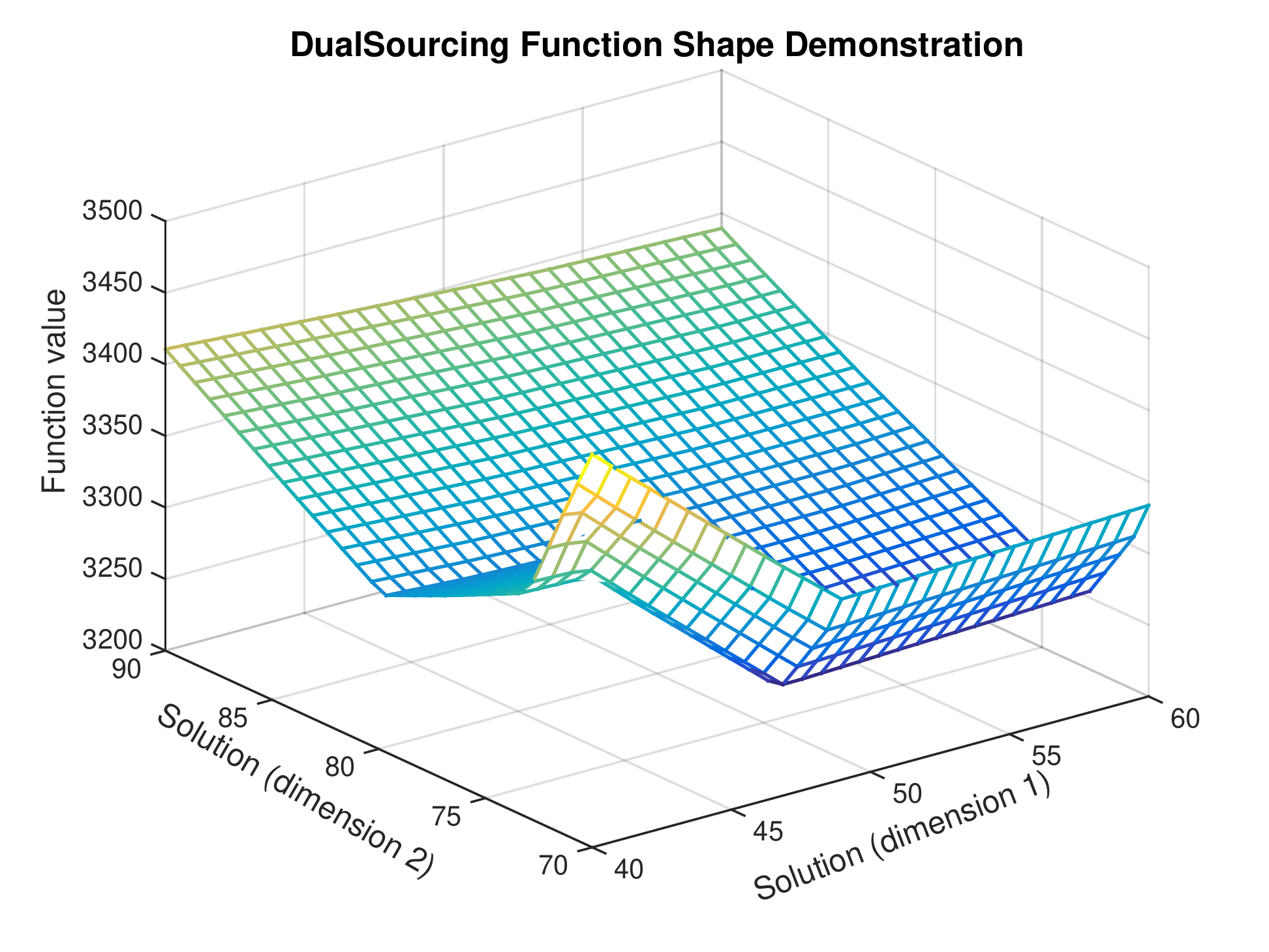}
					\end{subfigure}
	\caption{DualSourcing: (l) Average Performance w/ Bad Init., (r) Meshplot of the Objective Function}
	\label{fig:DualSourcing57}
\end{figure}

Our preliminary results for the problems EOQ and Dual Sourcing indicated that all of the algorithms had similar performance, quickly finding good solutions and then not improving further, as shown in Figure~\ref{fig:EOQ}.
We determined that the generated initial
solutions were already near optimal for these two problems, as seen in the vertical scales on
those plots.
To induce differences in the performances across algorithms, we reran the algorithms, intentionally generating poor initial solutions; see Figures~\ref{fig:EOQ10} and~\ref{fig:DualSourcing57}.

In Figure \ref{fig:EOQ10}, all algorithms except \GS perform better than \RS. \GS appears to fail due to the shape of the objective
function---the function slope is steep to the left of the
minimum but very flat to the right. Thus, when starting from a
relatively small initial solution, the algorithm can first take a
large step to a solution in the flat-slope area. Afterwards, \GS takes
very small steps back towards the optimal solution. The 
performance of \RS is highly dependent on the sampling distribution,
which in these examples is not well calibrated, leading to poor performance.

For the problem DualSourcing (Figure \ref{fig:DualSourcing57}), the extremely bad performance of \RS can again be explained by the poorly-claibrated distribution for generating solutions.
\SPSA performs particularly poorly.
The relatively weak performance of \NM on this problem is the result of a few macroreplications on which the algorithm failed to find the optimal solution.

\subsection*{MultiModal, POMDPController, Rosenbrock, Toll Network}
For these problems, most algorithms make rapid early progress and then stagnate.

\begin{figure}[H]
	\centering
	\begin{subfigure}[b]{0.48\textwidth}
		\includegraphics[width=\textwidth]{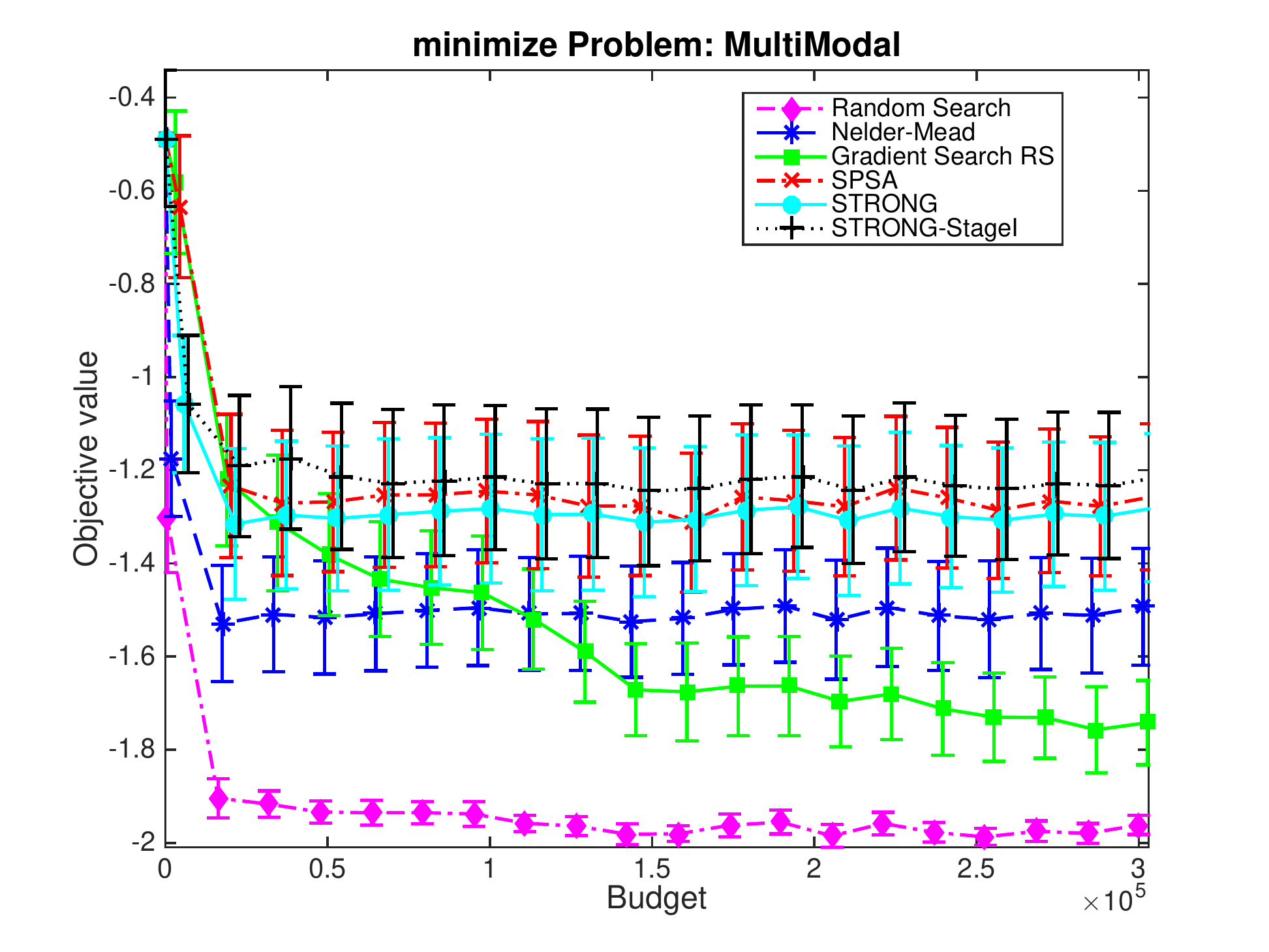}
					\end{subfigure}
	\begin{subfigure}[b]{0.48\textwidth}
		\includegraphics[width=\textwidth]{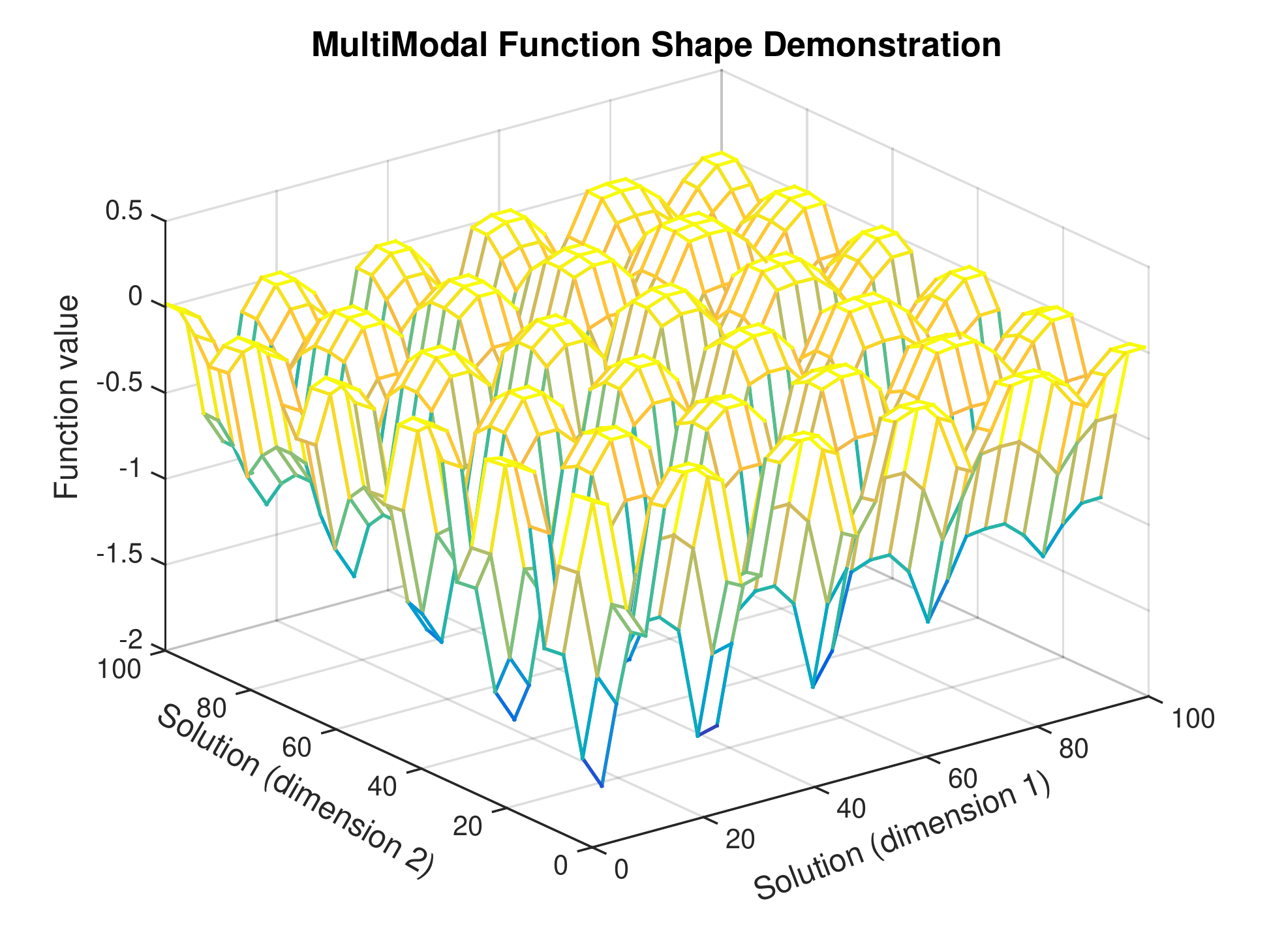}
					\end{subfigure}
	\caption{Multimodal: (l) Average Performance, (r) Meshplot of the Objective Function}
	\label{fig:MultiModal}
	\label{fig:MultiModal-Mesh}
\end{figure}

MultiModal is a 2-dimensional problem with 25 widely spaced local
optima (Figure~\ref{fig:MultiModal-Mesh}). Therefore, the local-search
algorithms quickly improve on the initial solution, but then fail to
make further progress. Meanwhile, the algorithms with restart, namely \RS and \GS, manage to identify better solutions. Compared to \RS, \GS is less
efficient because it only employs random restarts once it fails to
make progress, so it has fewer opportunities to restart. It only finds
the global optimum around half of the time, which helps to explain the
high variance of its performance in the results.

\begin{figure}[H]
	\centering
	\begin{subfigure}[b]{0.48\textwidth}
		\includegraphics[width=\textwidth]{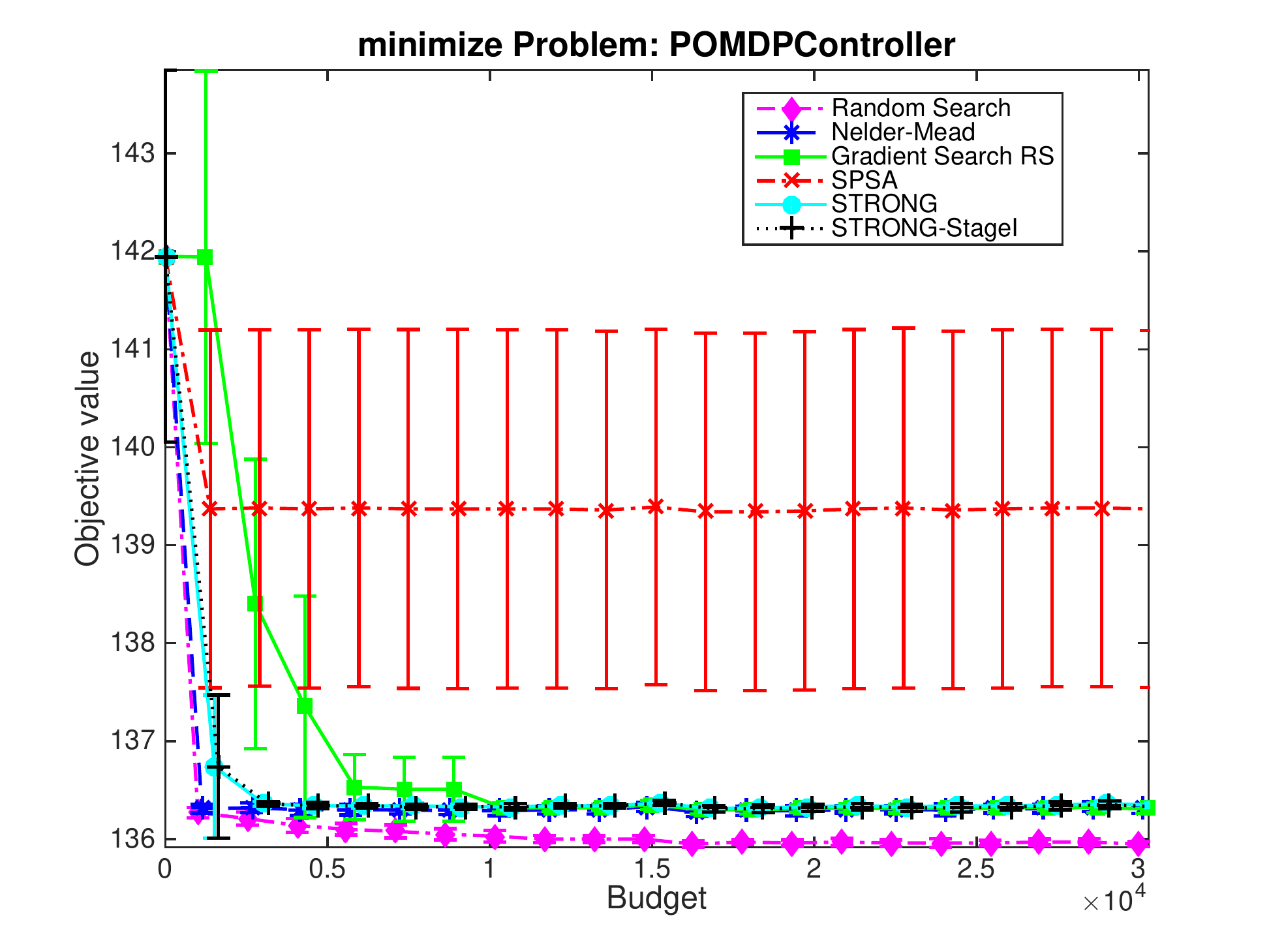}
					\end{subfigure}
	\begin{subfigure}[b]{0.48\textwidth}
		\includegraphics[width=\textwidth]{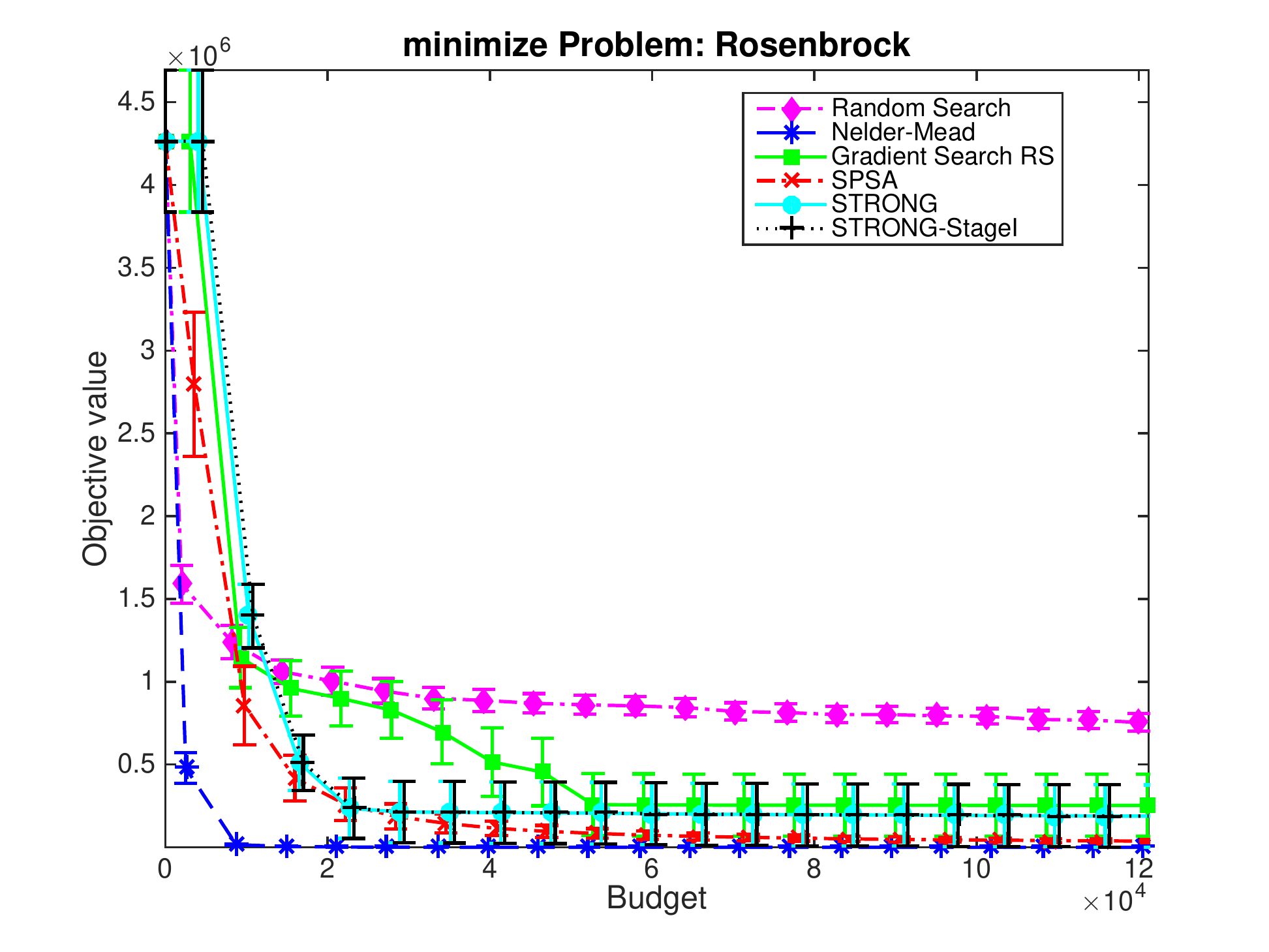}
					\end{subfigure}
	\caption{ (l) POMDPController: Average Performance, (r) Rosenbrock: Average Performance}
	\label{fig:POMDPController}
	\label{fig:Rosenbrock}
\end{figure}

We observed that the problem POMDPController was hard for most algorithms to solve (see Figure~\ref{fig:POMDPController}).
A plot of the objective function for the case $d=2$ shows that the objective function is
made up of two plateaus with a narrow valley running
through one of the plateaus. The optimal solution is at the bottom of
the valley. The valley is hard to find because it is such a small
region in the domain. In higher dimensions, such as the case $d=10$ that we tested, it is likely that good solutions are located in an even smaller region of the domain.
In addition, gradient-based methods struggle
because estimates of the gradient on the plateaus are close to zero
and therefore dominated by noise. This helps explain why \RS
outperforms the other methods. The mean performance of \SPSA is dominated by a few macroreplications on which it performed very poorly, perhaps because it produces very noisy gradient estimates that are based on
only two objective function evaluations.

The Rosenbrock problem has the highest dimension of all our problems. As expected, \RS performs poorly in searching a
high-dimensional space (Figure~\ref{fig:Rosenbrock}). \NM performs exceptionally well. Occasionally,
\GS reaches the boundary of the feasible region and the estimated
gradient indicates a search direction outside the boundary. We do not
use a projected-gradient algorithm, instead simply pinning solutions
back to the boundary, so \GS can stall on this problem. The two trust-region algorithms appear to struggle, perhaps for the same reason that
\GS struggles. As expected, \SPSA performs very well in this
high-dimensional problem, owing to its cheap gradient estimates.

\begin{figure} [H]
\centering
\includegraphics[width=0.48\textwidth]{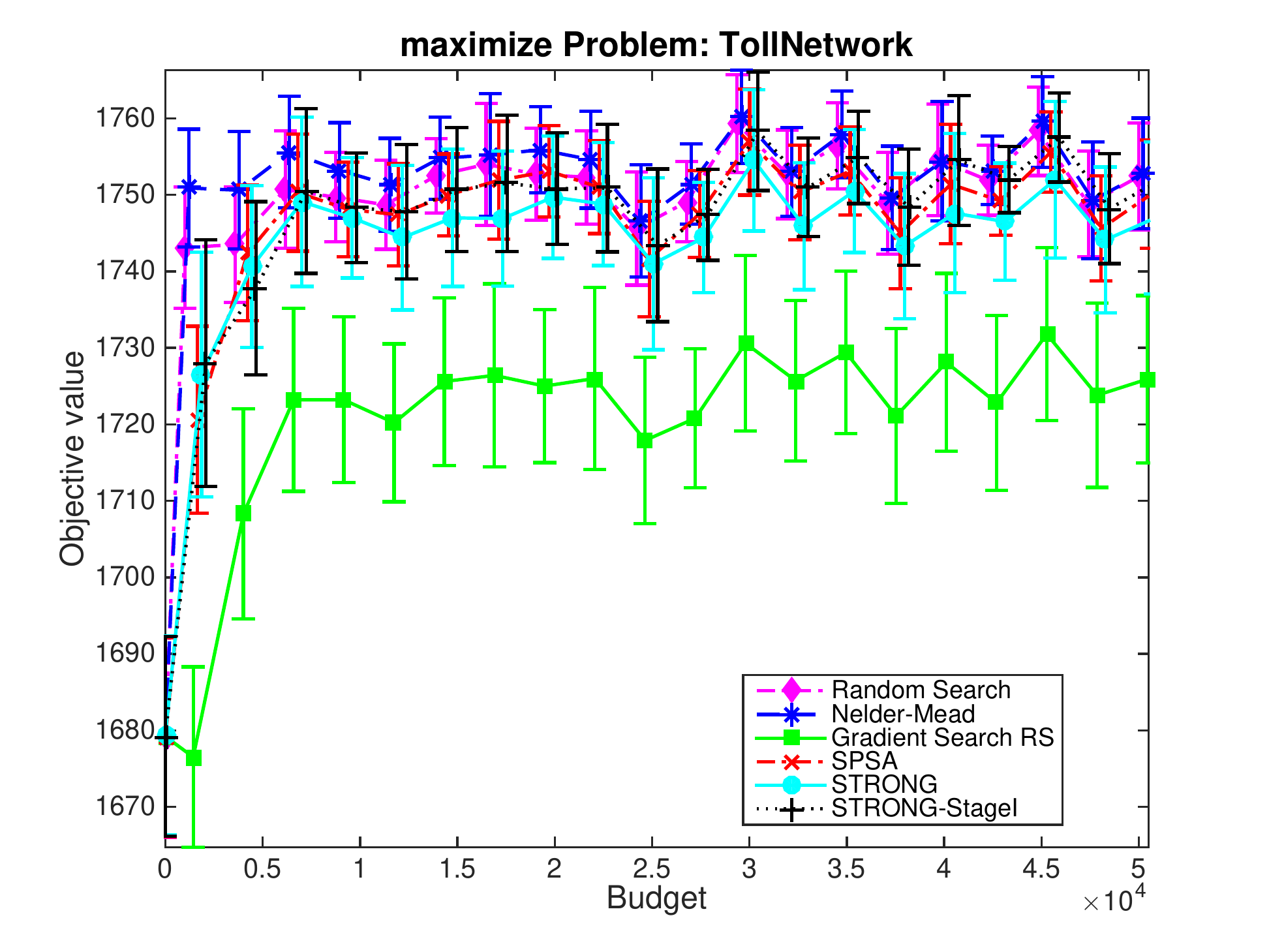}
\caption{TollNetwork: Average Performance}
\label{fig:TollNetwork}
\end{figure}

For the TollNetwork problem, \GS struggled relative to the other algorithms (Figure \ref{fig:TollNetwork}).
We do not have an explanation for why \GS performs so poorly on this problem.

\subsection*{Ambulance, FacilityLocation, ProductionLine, RoutePrices, SAN}
The algorithms have highly varied performance on these five problems. 

\begin{figure}[h]
	\centering
	\begin{subfigure}[b]{0.48\textwidth}
		\includegraphics[width=\textwidth]{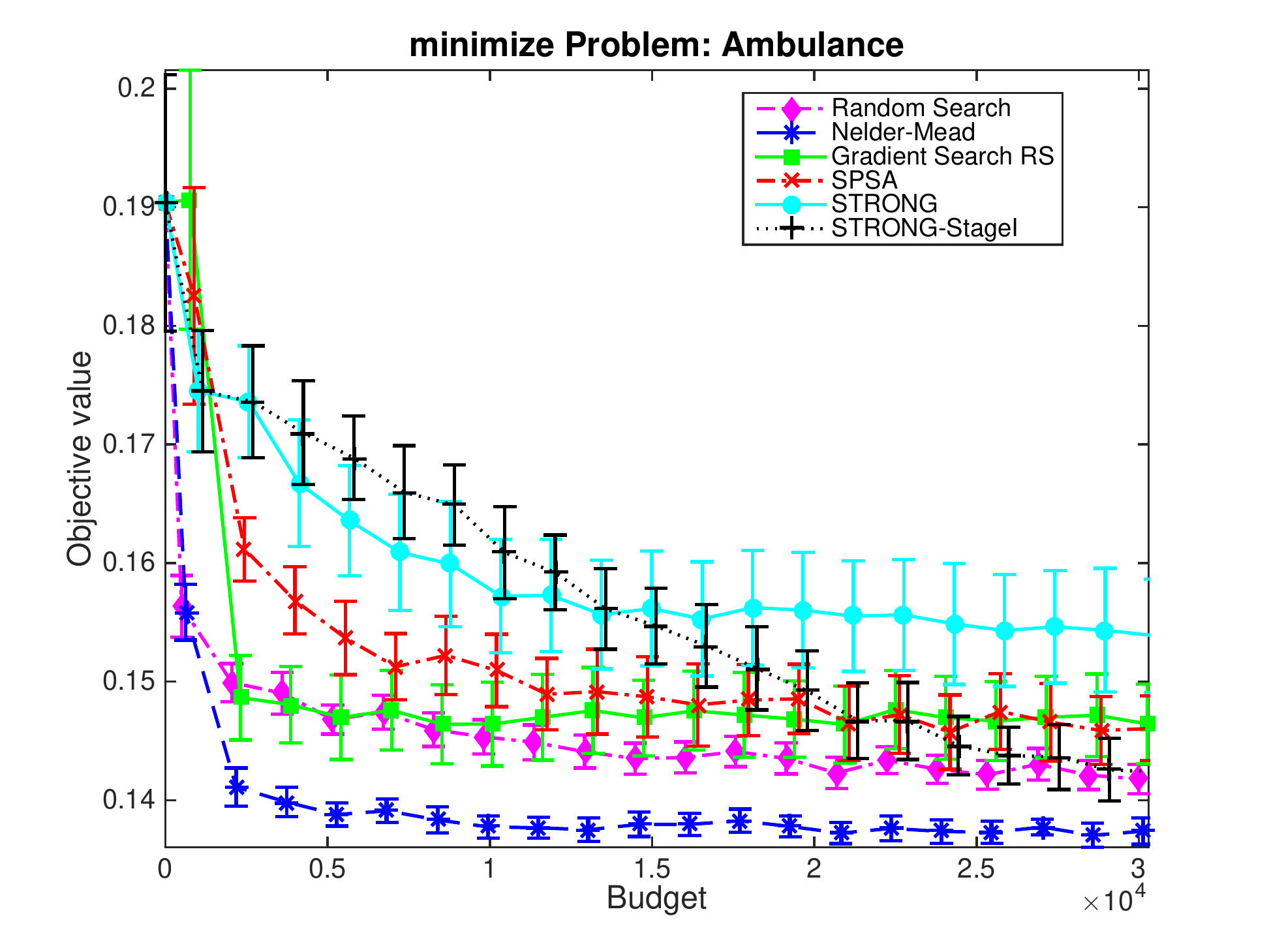}
					\end{subfigure}
	\begin{subfigure}[b]{0.48\textwidth}
		\includegraphics[width=\textwidth]{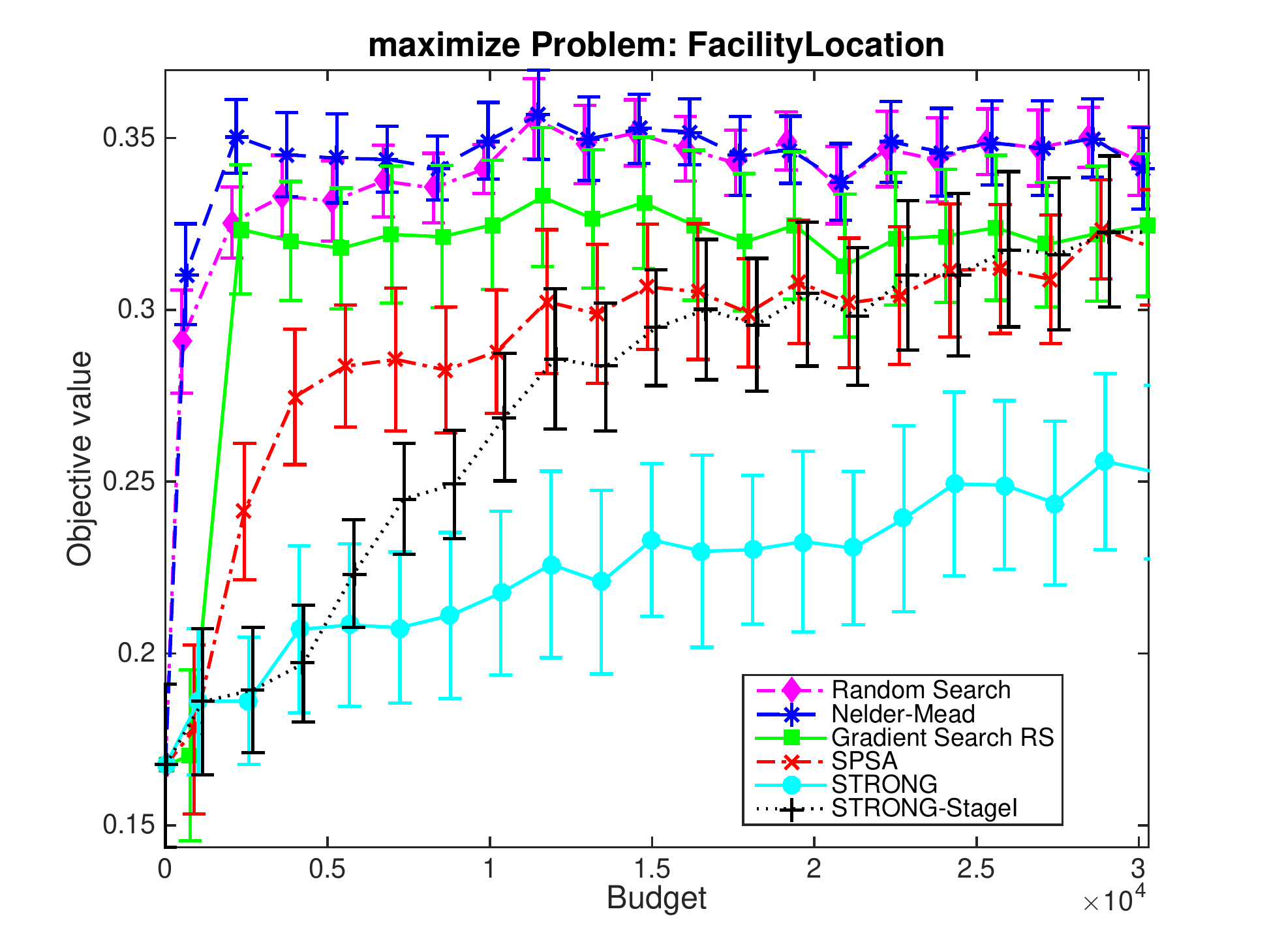}
					\end{subfigure}
	\caption{ (l) Ambulance: Average Performance, (r) FacilityLocation: Average Performance}
	\label{fig:Ambulance}
	\label{fig:FacilityLocation}
\end{figure}

The problems Ambulance and FacilityLocation have similar problem
structure, and the plots exhibit similar trends
(Figure~\ref{fig:Ambulance}). \NM is the stand-out performer, and
interestingly \RS performs almost as well. \STR probably struggles due
to the computational expense of iterations in Stage II.

\begin{figure}
	\centering
	\begin{subfigure}[b]{0.48\textwidth}
		\includegraphics[width=\textwidth]{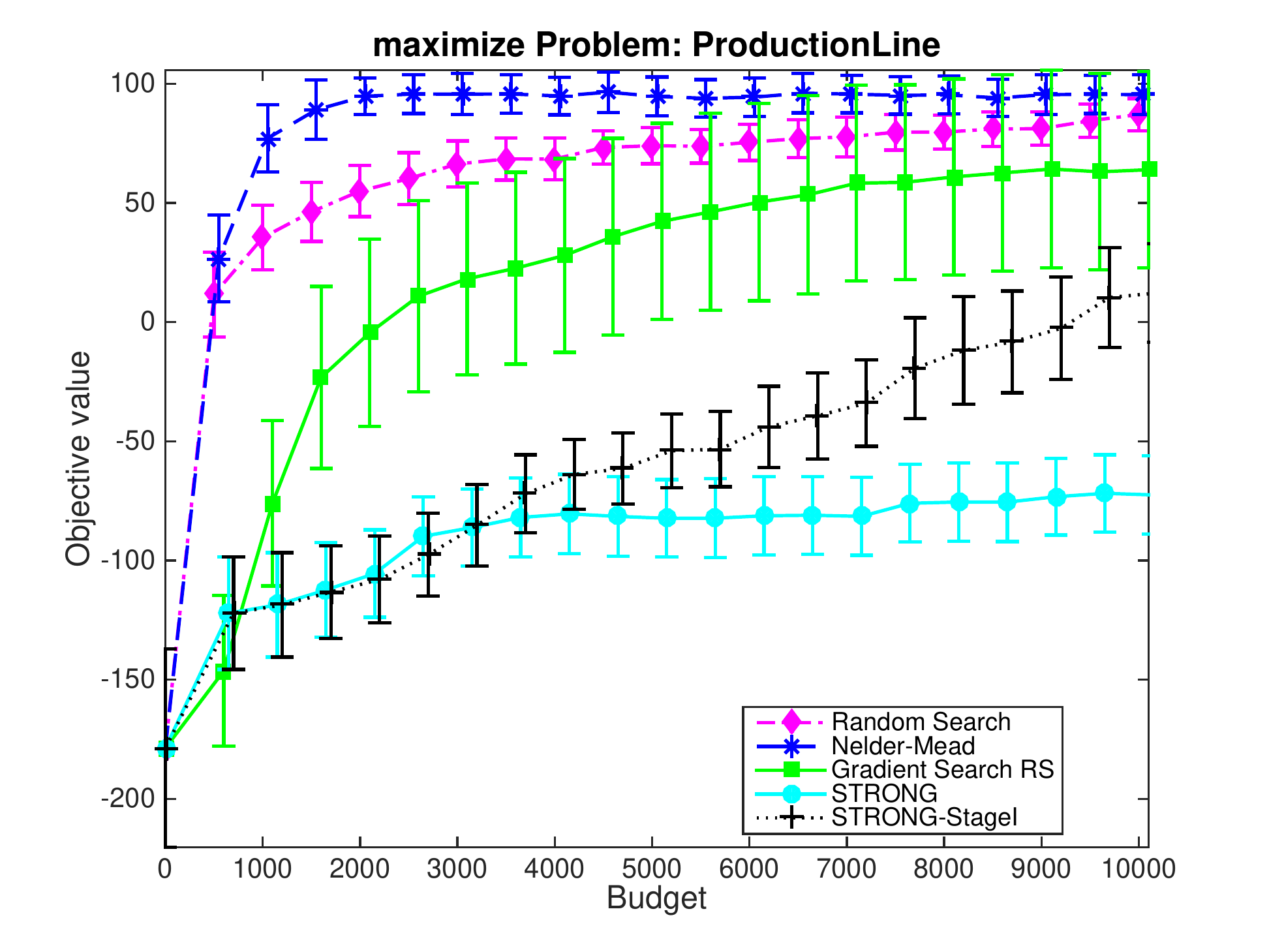}
					\end{subfigure}
	\begin{subfigure}[b]{0.48\textwidth}
		\includegraphics[width=\textwidth]{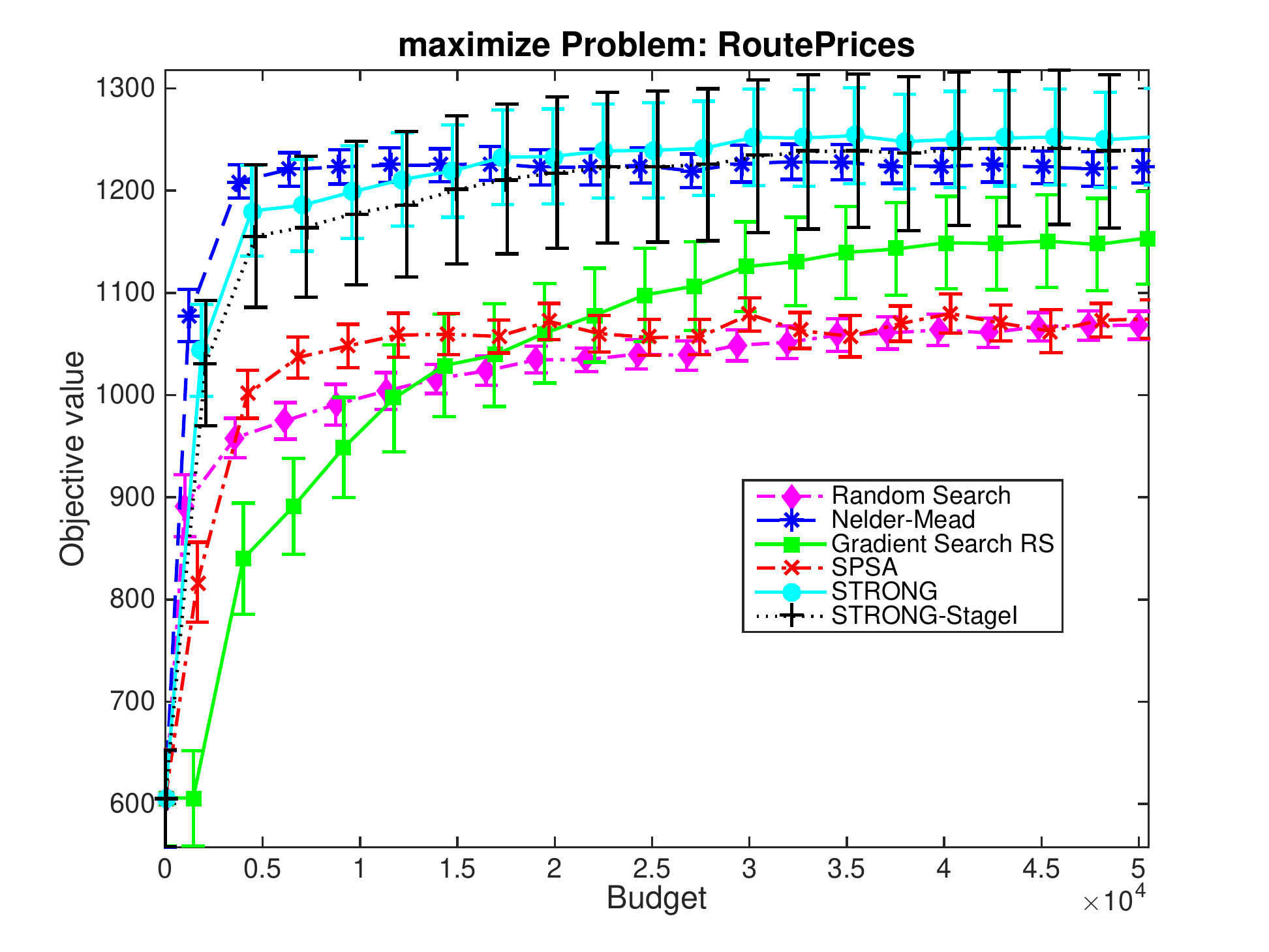}
					\end{subfigure}
	\caption{ (l) ProductionLine: Average Performance, (r) RoutePrices: Average Performance}
	\label{fig:ProductionLine}
	\label{fig:RoutePrices}
\end{figure}

For the problem ProductionLine, all algorithms perform much as they
did on the Ambulance problem, with the exception of \SPSA which is excluded from Figure
\ref{fig:ProductionLine} due to numerical problems in some macro
replications. The RoutePrices problem has dimension $d=12$, and
accordingly \RS performs poorly. \NM, \STR and \STRONE all perform
well on this problem. Of these three algorithms, \NM is the most
consistent performer across macroreplications as suggested
by the narrow confidence intervals, but both trust region methods
appear to identify a slightly better objective function for the
maximal budget. We do not have an
explanation for the poor performance of \GS and \SPSA on this problem.

\begin{figure} [H]
\centering
\includegraphics[width=0.48\textwidth]{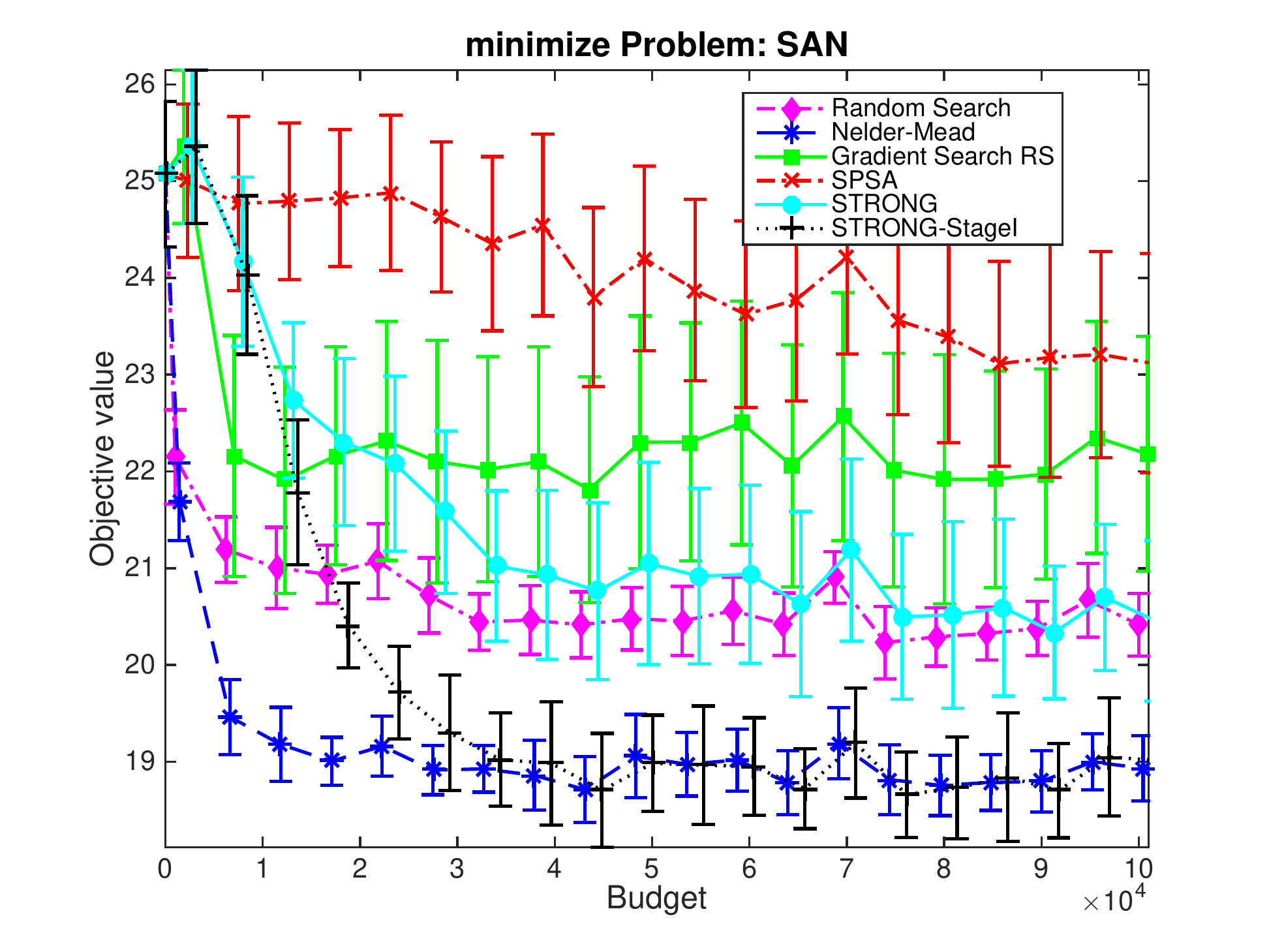}
\caption{SAN: Average Performance}
\label{fig:SAN}
\end{figure}

For the problem SAN (Figure~\ref{fig:SAN}), the strongest performers are \NM and
\STRONE. This problem is convex, so we would expect all algorithms to
do reasonably well. Nevertheless, \STR struggles, perhaps due to the
computational effort in Stage II. \GS also struggles, perhaps because
its test for when to perform a restart is too lenient, leading to many
unproductive restarts. \SPSA performs particularly poorly on this problem.

	\section{CONCLUSIONS}
	Perhaps the most important observation in our results is the robust
performance of \NM across almost all problems. It is highly deserving
of further study. Interestingly, \STRONE either compares similarly to,
or outperforms, \STR, calling into question the computationally
expensive Stage II of \STR. \RS performed better than expected, though
its performance suffers on higher-dimensional problems or on problems
where the sampling distribution is poorly calibrated. \SPSA performed
very well on our highest-dimensional problem (in 40 dimensions), but for most problems it
struggled relative to the other algorithms, including problems in
dimensions 10--13. Our implementation of \GS performed moderately
well, especially with its use of random restarts and may be worth
further development.

More generally, our sense is that low-dimensional problems are over-represented in the
SimOpt library. It is certainly conceivable that our observations above could
change for higher-dimensional problems.
The number of problems was small enough that we
could present detailed results for a single problem at a time, leading
to several insights. If we were to study many more problems, then this
presentation would be too cumbersome and a more succinct approach
would be necessary. To that end, further research on how to adapt
performance profiles to simulation optimization is needed.

In this study we only tackled continuous-variable problems that are
either unconstrained or box constrained (with upper and lower bounds
on the variables). 
In future research, we would like to compare algorithms designed for discrete or integer-ordered variables, like COMPASS \cite{hong:06} and R-SPLINE \cite{wang:13}.

\subsubsection*{ACKNOWLEDGMENTS}
This material is based upon work supported by the National Science Foundation under grant no.\ CMMI-1537394,
CMMI-1254298, CMMI-1536895, and IIS-1247696,
by the Air Force Office of Scientific Research under grant no.\ FA9550-12-1-0200, FA9550-15-1-0038, and FA9550-16-1-0046,
and by the Army Research Office under grant no.\ W911NF-17-1-0094.
%

		\bibliography{demobib}
		
\end{document}